\documentclass{amsart}
\usepackage[dvips]{graphicx}
\usepackage{amscd}
\usepackage{amsmath}
\usepackage{amsxtra}
\usepackage{amsfonts}
\usepackage{amssymb}
\newtheorem{theorem}{Theorem}[section]
\newtheorem{corollary}[theorem]{Corollary}
\newtheorem{lemma}[theorem]{Lemma}
\newtheorem{proposition}[theorem]{Proposition}
\theoremstyle{definition}
\newtheorem{definition}[theorem]{Definition}
\newtheorem{remark}[theorem]{Remark}

\theoremstyle{remark}

\renewcommand{\theclaim}{\textup{\theclaim}}

\newtheorem*{acknowledgements}{Acknowledgements}

\numberwithin{equation}{section}

\def\openone
{\mathchoice
{\hbox{\upshape \small1\kern-3.3pt\normalsize1}}
{\hbox{\upshape \small1\kern-3.3pt\normalsize1}}
{\hbox{\upshape \tiny1\kern-2.3pt\SMALL1}}
{\hbox{\upshape \Tiny1\kern-2pt\tiny1}}}
\makeatletter
\newbox\ipbox
\newcommand{\ip}[2]{\left\langle #1\mathrel{\mathchoice
{\setbox\ipbox=\hbox{$\displaystyle \left\langle\mathstrut #1#2\right\rangle$}
\vrule height\ht\ipbox width0.25pt depth\dp\ipbox}
{\setbox\ipbox=\hbox{$\textstyle \left\langle\mathstrut #1#2\right\rangle$}
\vrule height\ht\ipbox width0.25pt depth\dp\ipbox}
{\setbox\ipbox=\hbox{$\scriptstyle \left\langle\mathstrut #1#2\right\rangle$}
\vrule height\ht\ipbox width0.25pt depth\dp\ipbox}
{\setbox\ipbox=\hbox{$\scriptscriptstyle \left\langle\mathstrut #1#2\right\rangle$}
\vrule height\ht\ipbox width0.25pt depth\dp\ipbox}
} #2\right\rangle}
\newcommand{\diracb}[1]{\left\langle #1\mathrel{\mathchoice
{\setbox\ipbox=\hbox{$\displaystyle \left\langle\mathstrut #1\right.$}
\vrule height\ht\ipbox width0.25pt depth\dp\ipbox}
{\setbox\ipbox=\hbox{$\textstyle \left\langle\mathstrut #1\right.$}
\vrule height\ht\ipbox width0.25pt depth\dp\ipbox}
{\setbox\ipbox=\hbox{$\scriptstyle \left\langle\mathstrut #1\right.$}
\vrule height\ht\ipbox width0.25pt depth\dp\ipbox}
{\setbox\ipbox=\hbox{$\scriptscriptstyle \left\langle\mathstrut #1\right.$}
\vrule height\ht\ipbox width0.25pt depth\dp\ipbox}
}\right. }
\newcommand{\dirack}[1]{\left. \mathrel{\mathchoice
{\setbox\ipbox=\hbox{$\displaystyle \left.\mathstrut #1\right\rangle$}
\vrule height\ht\ipbox width0.25pt depth\dp\ipbox}
{\setbox\ipbox=\hbox{$\textstyle \left.\mathstrut #1\right\rangle$}
\vrule height\ht\ipbox width0.25pt depth\dp\ipbox}
{\setbox\ipbox=\hbox{$\scriptstyle \left.\mathstrut #1\right\rangle$}
\vrule height\ht\ipbox width0.25pt depth\dp\ipbox}
{\setbox\ipbox=\hbox{$\scriptscriptstyle \left.\mathstrut #1\right\rangle$}
\vrule height\ht\ipbox width0.25pt depth\dp\ipbox}
} #1\right\rangle}

\input cyracc.def

\begin{document}
\title[Inner Amenability on the Generalized Thompson Group]{The Inner Amenability of the Generalized Thompson Group }
\author{Gabriel Picioroaga}
\address{Department of Mathematics\\
The University of Iowa\\
14 MacLean Hall\\
Iowa City, IA 52242-1419\\
U.S.A.}
\email{gpicioro@math.uiowa.edu}
\thanks{}

\begin{abstract} In this paper we prove that the general version, $F(N)$ of the Thompson group is 
inner amenable. As a consequence we generalize a result of P.Jolissaint. To do so, we prove first that 
$F(N)$ together with a normal subgroup are i.c.c (infinite conjugacy classes) groups. Then, we 
investigate the relative McDuff property out of which we extract property $\Gamma$ for the group von 
Neumann algebras involved. By a result of E.G.Effros, $F(N)$ follows inner amenable. 
\end{abstract}\maketitle

\section{Introduction}
The Thompson group $F$ can be regarded as the group of piecewise-linear,\\ orientation-preserving 
homeomorphisms of the unit interval which have breakpoints only at dyadic points and on intervals of 
differentiability the slopes are powers of two. The group was discovered in the '60s by Richard Thompson 
and in connection with the now celebrated groups $T$ and $V$ it led to the first example of a finitely 
presented infinite simple group. Since then, these groups have received considerable applications in such 
fields as homotopy theory or operator algebras. The group $V$ has been generalized by Higman (\cite{Hig}) and $F$ by 
Brown (\cite{Br}) and M. Stein (\cite{St}). The generalized Thompson group we study in this paper 
corresponds to $F(p)$ in \cite{St} (or $F_{p,\infty}$ in \cite{Br}). 
\par In 1979 Geoghegan conjectured that $F$ is not amenable. This problem is still open and of great 
importance for group theory: either outcome will produce a strict inclusion related to groups satisfying 
certain properties (i.e. there would exist a finitely presented group that is not amenable and does not 
contain a free subgroup of rank 2 or else there would be a finitely presented amenable group which is 
not elementary amenable). 
\par Inner amenability (a larger property than amenability eventhough the definition is just a slightly 
variation of the amenability one) was introduced by Effros. He also observed that if the group von 
Neumann algebra corresponding to a group $G$ is a $II_1$ factor having the property $\Gamma$ of Murray 
and von Neumann then $G$ is inner amenable. In \cite{Jol1} Jolissaint proved that the Thompson group is 
inner amenable, then in \cite{Jol} he proved more: the $II_1$ factor associated with the Thompson group 
has the so-called relative McDuff property; in particular property $\Gamma$ is satisfied.\\
In this paper we prove that for any integer $N\geq 2$ the generalized Thompson group $F(N)$ satisfies 
the relative McDuff property and therefore it is inner amenable, hence the title. To prove the McDuff 
property we will use two results from \cite{Jol}. We remark that inner amenability could be 
proved directly using the ideas in \cite{Jol1} and the infinite presentation of $F(N)$, however relative 
McDuff property is a stronger result.
In the next section we prepare some basics on the Thompson groups and group von Neumann algebras. In the 
last section we prove the main result of the paper and its corollary. We conclude with a question 
whose answer might connect the theory of $II_1$ factors and the (non)amenability of $F$.

\section{Background}

\begin{definition}For $N\in\mathbb{N}$, $N\geq 2$, the Thompson group $F(N)$ is the set of piecewise 
linear homeomorphisms from the closed unit interval $[0,1]$ to itself that are differentiable except at 
finitely many $N$-adic rationals and such that on intervals of differentiability the derivatives are 
powers of $N$. 
\end{definition}
\par
In \cite{Br} finite and infinite presentations of $F(N)$ are given. For example
$$F(N)=\left< x_0, x_1, ...x_i,...|\mbox{ }x_jx_i=x_ix_{j+N-1}\mbox{, }i<j\mbox{ }\right>$$

We will not make use of them here, our arguments being based on some special elements of $F(N)$ 
(see \cite{DutP}) 
$$A_{d,p}(x)=\left\{
\begin{array}{lr}
x/N^p,&\mbox{ } 0\leq x\leq d\\
x-d+d/N^p,&\mbox{ } d\leq x\leq 1-d/N^p\\
N^px+1-N^p,&\mbox{ } 1-d/N^p\leq x\leq 1
\end{array}\right. $$
where $d$ is a N-adic, $p\in\mathbb{Z}$ such that $d/N^p<1$. \\
Next, we introduce two subgroups of $F(N)$. 
Let $$F^{'}:=\{f\in F(N)\mbox{ }|\mbox{ }f_{|[0,\epsilon]}=id\mbox{, }f_{|[\delta,1]}=id,\mbox{ }0<\epsilon\mbox{,}\delta<1\}$$
 and 
the intermediate subgroup
 $$D:=\{f\in F(N)\mbox{ }|\mbox{ }f_{|[\delta,1]}=id,\mbox{ }0<\delta<1\}$$
These, of course are the same ones considered in \cite{Jol} for $N=2$. It is not hard to see that 
$F^{'}$ and $D$ are normal subgroups of $F(N)$ . Actually, when $N=2$, $F^{'}$ is the commutator 
subgroup (see \cite{Can} for $N=2$ and \cite{Br} for the general case). \\

\par A von Neumann algebra can be thought of as a * subalgebra of bounded operators on some 
fixed (separable) Hilbert space, that is closed with respect to the weak topology.  As shown by John von Neumann 
the building blocks of the theory are the so-called factors. 
A von Neumann algebra $M$ inside $B(H)=$ the space of all bounded operators on some Hilbert space $H$ is 
called a factor if $M\cap M^{'}=\mathbb{C}$, where $M^{'}$ represents the set of all bounded operators 
commuting with any element of $M$. There are three types of factors:\\ 
type I, when $M$ admits minimal projections, e.g. $n\times n$ matrices over $\mathbb{C}$ or the whole $B(H)$;\\
type II, no minimal projection and there exists a unique (semi) finite trace. The factor is called $II_1$ 
when the trace takes only finite values;\\
type III means not of type I or II.\\
For more on factors and von Neumann algebras we refer the reader to the book \cite{StZ}.  
\par Since the beginning of the theory there has been a fruitful interplay between group theory and factors. 
For example, if $G$ is a countable discrete group with infinite conjugacy classes (i.c.c.) then the left 
regular representation of $G$ on $l^2(G)$ gives rise to a $II_1$ factor, the group von Neumann algebra $\mathcal{L}(G)$, 
as follows:\\
Let $l^2(G)=\left\{ \psi:G\rightarrow\mathbb{C}\mbox{ }|\mbox{}\sum_{g\in G}|\psi(g)|^2<\infty\right\}$ 
endowed with the scalar product 
$$\left<\phi,\psi\right>:=\sum_{g\in G}\phi(g)\overline{\psi(g)}$$
Notice that the Hilbert space $l^2(G)$ is generated by the countable colection of vectors 
$\left\{\delta_g\mbox{ }|\mbox{}g\in G\right\}$. Also, an element $g\in G$ defines a unitary operator $L_g$, 
on $l^2(G)$  as follows: 
$L_g(\psi)(h)=\psi(g^{-1}h)$, for any $\psi\in l^2(G)$ and any $h\in G$. (Sometimes, to not burden 
the notation we will write just $g$ instead of $L_g$). Now, $\mathcal{L}(G)$, 
the von Neumann algebra generated by $G$ is obtained by taking the wo-closure in $B(l^2(G))$ of the linear span of the set 
$\left\{L_g \mbox{ }|\mbox{}g\in G\right\}$. It is a routine exercise to show that $\mathcal{L}(G)$ is a 
factor provided any element of $G$ has infinite conjugacy classes, and it is of type $II_1$. 
The map defined  by 
$\mbox{tr}(L)=\left<L(\delta_e),\delta_e\right>$, where $e\in G$ is the neutral element and $L\in\mathcal{L}(G)$  
is a faithful, normal trace. The canonical trace also determines the Hilbertian norm $|x|_2=\mbox{tr}(x^{*}x)^{1/2}$. 
It is easy to see that $\mbox{tr}(g)=0$ for $g\neq e$ and $\mbox{tr}(e)=1$. Also, for $g\neq h$, $|g-h|_2=2^{1/2}$.\\
\\
The following is an equivalent definition of (inner)amenability:
\begin{definition}Let $G$ be a (countable discrete) group. If there exists a mean $f$ on 
the algebra $l^{\infty}(G-\{e\})$, invariant under the action $$(gf)(h)=f(g^{-1}h)$$
then $G$ is amenable. If the action is taken with respect to conjugation then $G$ is inner amenable. 
\end{definition}
Let $M$ be a $II_1$ factor and $\mbox{tr}$ its normal, faithful trace. 
\begin{definition}i) A sequence $(x_n)_n\in l^{\infty}(\mathbb{N},M)$ is a central sequence if 
$$\lim_{n\rightarrow\infty}|[x,x_n]|_2=0$$
for any $x\in M$, where $[x,y]=xy-yx$.\\
ii) Two central sequences $(x_n)$ and $(y_n)$ in $M$ are equivalent if 
$$\lim_{n\rightarrow\infty}|x_n-y_n|_2=0$$  
iii)A central sequence is trivial if it is equivalent to a scalar sequence.\\
iv)$M$ has property $\Gamma$ of Murray and von Neumann if there exists in $M$ a non trivial central sequence.
\end{definition}
In \cite{Efr} it is shown that if the $II_1$ factor $\mathcal{L}(G)$ satisfies property $\Gamma$ 
then $G$ is inner amenable. A stronger property than $\Gamma$ is the relative McDuff property: a factor 
$M$ is McDuff if $M$ is isomorphic to $M\otimes R$, where $R$ is the hyperfinite $II_1$ factor ($R$ can be 
viewed as $\mathcal{L}(G)$ for an amenable, countable, discrete group $G$). D.Bisch extended this property 
to pairs of $II_1$ factors $1\in N\subset M$:
\begin{definition} The pair $N\subset M$ has the relative McDuff property if there exists an isomorphism
$\Phi:M\rightarrow M\otimes R$ such that $\Phi(N)=N\otimes R$.
\end{definition}
\par Finally, we prepare the two results from \cite{Jol} that we are going to use. First, notice that an 
inclusion $H\subset G$ of groups determines naturally an embedding 
$1\in\mathcal{L}(H)\subset\mathcal{L}(G)$ of factors. Also, a semidirect product $H\rtimes_{\alpha}G$ 
of groups translates into a crossed product  $\mathcal{L}(H)\rtimes_{\alpha}G$ in the 
realm of factors (for more on crossed products see for example \cite{vDa}). 
\begin{proposition}( Proposition 2.4. in \cite{Jol}) Let $G$ a countable i.c.c. group and let $H$ be an 
i.c.c. subgroup of $G$ with the following property: for every finite subset $E$ of $G$ there exist 
elements $g$ and $h$ in $H-\{e\}$ such that\\
(1) $xg=gx$ and $xh=hx$ for every $x\in E$;\\
(2) $gh\neq hg$.\\
Then the pair $\mathcal{L}(H)\subset\mathcal{L}(G)$  has the relative McDuff property. 
\end{proposition}
\begin{definition}Let $M$ be a type $II_1$ factor and let $\theta$ be an automorphism of $M$.\\
i) $\theta$ is centrally trivial if one has for every central sequence $(a_n)_n$ 
$$\lim_{n\rightarrow\infty}|\theta(a_n)-a_n|_2=0$$
ii)If $G$ is a countable group and if $\alpha$ is an action of $G$ on $M$, then $\alpha$ is called 
centrally free if $\alpha_g$ is not centrally trivial for every $g\in G-\{e\}$.
\end{definition}
\begin{proposition}(Proposition 2.6. in \cite{Jol}) Let $N$ be a McDuff factor of type $II_1$ with 
separable predual, let $G$ be an amenable countable group and let $\alpha$ be a centrally free action 
of $G$ on $N$. Then the pair $N\subset N\rtimes_{\alpha}G$ has the relative McDuff property.
\end{proposition}
\section{Main result}

As in \cite{Jol} we will establish that the  pairs $\mathcal{L}(F^{'})\subset\mathcal{L}(D)$ and 
$\mathcal{L}(D)\subset\mathcal{L}(F(N))$ have the relative McDuff property. In doing so we will first 
make sure that the group von Neumann algebras involved are $II_1$ factors. 

\begin{theorem}\label{t1} $F^{'}$ and $F(N)$ are both i.c.c. groups, therefore the von Neumann algebras 
$\mathcal{L}(F^{'})$ and $\mathcal{L}(F(N))$ are $II_1$ factors.
\end{theorem}
\begin{proof}The nice argument that $F^{'}$ is i.c.c. is due to Dorin Dutkay: notice first that for any 
non-trivial $g\in F^{'}$ there is a unique $\epsilon$ such that $g_{[0,\epsilon]}=\mbox{id}$ and for any 
neighborhood $V$ of $\epsilon$ there is $x\in V$, $x>\epsilon$, $g(x)\neq x$. We call this unique value $\epsilon_g$. Now, we prove that for any 
$h\in F$ 
\begin{equation}\label{e1}
\epsilon_{hgh^{-1}}=h(\epsilon_g)
\end{equation}
$\epsilon_{hgh^{-1}}$ makes sense because $g\in F^{'}$ implies $hgh^{-1}\in F^{'}$. 
Let $x\in [0, h(\epsilon_g)]$. Then $gh^{-1}(x)=h^{-1}(x)$, hence $hgh^{-1}(x)=x$. The maximality of 
$\epsilon_{hgh^{-1}}$ implies 
$$\epsilon_{hgh^{-1}}\geq h(\epsilon_g)$$
For the reversed inequality apply the above with the substitutions $h \rightarrow h^{-1}$ and 
$g\rightarrow hgh^{-1}$. Therefore 
(\ref{e1}) holds. But, for a fixed non-trivial $g$,  the set $\{h(\epsilon_g)\mbox{ }|\mbox{ }h\in F^{'}\}$ 
is infinite, so that by (\ref{e1}) above, the conjugacy class of $g$ is infinite. In conclusion,  $F^{'}$ 
is i.c.c.
\par To prove $F(N)$ is i.c.c we are going to use the elements $A_{d,p}$. Notice first 
$$A^{-1}_{d,p}(x)=\left\{
\begin{array}{lr}
xN^p,&\mbox{ } 0\leq x\leq d/N^p\\
x+d-d/N^p,&\mbox{ } d/N^p\leq x\leq 1-d\\
(x+N^p-1)/N^p,&\mbox{ } 1-d\leq x\leq 1
\end{array}\right.$$
We assume $f\in F(N)$ is non-trivial and that in a neighborhood of $x=0$ its slope is positive. Once we 
prove that its conjugacy class is infinite the class of its inverse follows infinite, so there is no loss 
of generality in assuming a positive slope around the origin. For 
any $p\in\mathbb{Z}$, $p\neq 0$ we will find a large $\alpha$ such that for any $k>l>\alpha$: 
\begin{equation}\label{e2}
A_{d,p}^{-1}fA_{d,p}\neq A_{\tilde{d},p}^{-1}fA_{\tilde{d},p}
\end{equation}
where $d=1/N^k$ and $\tilde{d}=1/N^l$.
This relation clearly shows that the conjugacy class of $f$ is infinite. \\
From the definition of $F(N)$ there is a $N$-adic $d_1$ such that $f_{|[0,d_1]}(x)=N^nx$, where, by our 
assumption $n>0$ (the case $n=0$ is taken care of by the previous argument, i.e. when $f$ is 
trivial around $x=0$).  Choose now $\alpha$ large enough such that the following inequalities hold:  
$$\frac{1}{N^{\alpha+p}}<d_1\mbox{,  }\frac{N^n}{N^{\alpha+p}}<1-\frac{1}{N^{\alpha}}$$
For $k>l>\alpha$ consider the N-adic numbers $d$ and $\tilde{d}$. If (\ref{e2}) were not true then 
evaluating at $x=d=1/N^k$ we obtain   
$$A^{-1}_{d,p}f(\frac{1}{N^{k+p}})=A^{-1}_{\tilde{d},p}f(\frac{1}{N^{k+p}})$$
Because $1/N^{k+p}\in [0,d_1]$ the equality becomes
\begin{equation}\label{e3}
A^{-1}_{d,p}(\frac{N^n}{N^{k+p}})=A^{-1}_{\tilde{d},p}(\frac{N^n}{N^{k+p}})
\end{equation}
Because of the choices of $n$, $\alpha$, $k$ and $l$ we have
$$\frac{1}{N^{i+p}}<\frac{N^n}{N^{i+p}}<\frac{N^n}{N^{\alpha+p}}<1-\frac{1}{N^{\alpha}}<1-\frac{1}{N^i}$$
where $i\in\{k,l\}$.
This shows that 
$$x:=\frac{N^n}{N^{k+p}}\in [\frac{d}{N^p},1-d]$$
Using formula of $A^{-1}_{d,p}$, equation (\ref{e3}) can be rewritten
\begin{equation}\label{e4}
x+d-\frac{d}{N^p}=A^{-1}_{\tilde{d},p}(x)
\end{equation}
Because of the way $\alpha$ has been chosen we get $x<1-\tilde{d}$, so that there are only two cases to discuss: \\
$x\in [\frac{\tilde{d}}{N^p},1-\tilde{d}]$. Relation (\ref{e4}) easily implies $d=\tilde{d}$, which of course is 
not allowed.\\
$x\in [0, \frac{\tilde{d}}{N^p}]$. Using formula for $A^{-1}_{\tilde{d},p}$ on $[0, \frac{\tilde{d}}{N^p}]$ 
and putting $x=\frac{N^n}{N^{k+p}}$ and $d=\frac{1}{N^k}$ back in (\ref{e4}) we obtain
$$\frac{N^n}{N^{k+p}}+\frac{1}{N^k}-\frac{1}{N^{k+p}}=\frac{N^n}{N^k}$$
which reduces (using $n\neq 0$) to $N^{k+p}=N^k$.  This would imply $p=0$, a value that we avoid. \\
In conclusion, (\ref{e2}) is true and we finish the proof.    
\end{proof}
We prove the following lemma, useful for scaling down graphs of elements in $F(N)$ and still remaining in $F(N)$:
\begin{lemma}Let $0<\delta<\epsilon<1$, $N$-adic numbers  such that $\epsilon-\delta\in N^{\mathbb{Z}}$. 
Then there exists $f\in F(N)$ with $f_{|[0,\delta]}(x)=x$ and $f_{|[\epsilon,1]}(x)=x$ and $f$ has no fixed 
points in $(\delta,\epsilon)$. 
\end{lemma}
\begin{proof}Take $r:[\delta, \epsilon]\rightarrow [0,1]$ defined by $r(x)=\frac{x}{\epsilon-\delta}-\frac{\delta}{\epsilon-\delta}$. 
For $d$, a non-zero $N$-adic and $p\neq 0$, consider the following homeomorfism 
$$f(x)=\left\{
\begin{array}{lr}
x,&\mbox{ } 0\leq x\leq\delta\\
r^{-1}A_{d,p}r(x),&\mbox{ } \delta\leq x\leq\epsilon\\
x,&\mbox{ } \epsilon\leq x\leq 1
\end{array}\right.$$
When the derivative exists, $f^{'}(x)=A_{d,p}^{'}(r(x))$, for $x\in (\delta,\epsilon)$. Thus we obtain 
$f\in F(N)$. Also, there can be no fixed point of $f$ in $(\delta,\epsilon)$ as $A_{d,p}$ has no other fixed points besides $0$ and $1$. 
\end{proof}
Next, we are going to check the hypotheses of Proposition 2.4 from \cite{Jol}: 
\begin{proposition}\label{pr1}For every finite subset $\{g_1,g_2,...,g_n\}$ of $D$ there exist non-trivial distinct 
elements $g$ and $h$ of $F^{'}$ such that \\
(1) $g_ig=gg_i$ and $g_ih=hg_i$ for all $i\in\{1,2,...n\}$ \\
(2) $hg\neq gh$.
\end{proposition}
\begin{proof} Because $g_i\in D$, there exists $\delta$ such that $g_{i|[\delta,1]}=\mbox{id}$ for all 
$i\in\{1,2,...,n\}$. Also, we may take $\delta$ $N$-adic. For $\epsilon_1$ chosen such that 
$\epsilon_1-\delta\in N^{\mathbb{Z}}$ apply the Lemma above: there exists $g\in F^{'}$ with $g_{[0,\delta]}=\mbox{id}$, 
having no fixed points inside $(\delta,\epsilon_1)$.  An easy check shows $g_ig=gg_i$. Now, 
for $\epsilon_2>\epsilon_1$ such that $\epsilon_2-\delta\in N^{\mathbb{Z}}$ we find, using Lemma again, 
a $h\in F^{'}$ with the very same properties. Therefore (1) is satisfied. Now, for any $f\in F^{'}$ let 
$\epsilon_h$ the smallest $\epsilon$ such that $f_{|[\epsilon,1]}=$id and for any neighborhood $V$ of 
$\epsilon$ there exists a $x\in V$, $x<\epsilon$ with $f(x)\neq x$. For $g$ and $f$ found above we clearly 
have $\epsilon_1=\epsilon_g$ and $\epsilon_2=\epsilon_h$. Moreover, as in the proof of Theorem \ref{t1} 
the following equality holds true: $\epsilon_{hgh^{-1}}=h(\epsilon_g)$. If (2) were not satisfied then 
we would obtain $\epsilon_1=h(\epsilon_1)$, which contradicts the fact that $h$ has no fixed point 
inside $(\delta,\epsilon_2)$. 
\end{proof}
\begin{corollary}The pair  $\mathcal{L}(F^{'})\subset\mathcal{L}(D)$ has the relative McDuff property.
\end{corollary}
\begin{proof} Apply Proposition 2.4. of \cite{Jol}  
\end{proof} 
\begin{remark} We will continue on the ideas in \cite{Jol} to show that 
the pair \\ $\mathcal{L}(D)\subset\mathcal{L}(F(N))$ has the relative McDuff property. Notice first 
that the general Thompson group can be realized as a semiproduct $D\rtimes_{\alpha}\mathbb{Z}$ where 
the action $\alpha$ is defined as follows: choose $x_0\in F(N)$ such that its first piece of graph 
is trivial and the slope of the last piece is $N$ (it is elementary to construct such an element in 
$F(N)$, see the proof of Proposition \ref{pr2} below). Then the action $\alpha(n)(f)=x_0^nfx_0^{-n}$ is 
well defined on $D$. Also, any element of $F(N)$ can be written as $fx_0^{-n}$ for some $f\in D$ and 
$n\in\mathbb{Z}$, therefore the map $(f,n)\in D\rtimes_{\alpha}\mathbb{Z}\rightarrow fx_0^{-n}$ is a 
group isomorphism. For a central sequence in $\mathcal{L}(D)$ we choose the unitary operators 
corresponding to the following sequence $(a_n)_n\subset D$: let $(d_n)_n$ and $(\overline{d_n})_n$ two 
sequences of $N$-adic numbers in $[0,1]$ such that $d_n<\overline{d_n}$ are consecutives with $d_n\rightarrow 1$. 
Applying the scalling-down lemma we obtain a non-trivial $a_n\in F^{'}$. Moreover, for any $g$ in a finite 
subset of $D$ and for large $n$, $a_n$ commutes with $g$ (see the proof of Proposition \ref{pr1}).\\
Notice that $F^{'}$ is inner amenable: in particular, by the Corolarry above, $\mathcal{L}(F^{'})$ has 
property $\Gamma$, so that me may apply the result in \cite{Efr} to conclude $F^{'}$ is inner amenable. 
\end{remark}
\begin{theorem}The pair $\mathcal{L}(D)\subset\mathcal{L}(F(N))$ has the relative McDuff property. 
\end{theorem}
\begin{proof}We will make use of Proposition 2.6 in \cite{Jol}: it suffices to check that the action 
$\alpha$ is centrally free. Having already a central sequence given by $(a_n)_n$, it is enough to 
show $\mbox{lim}_n|\alpha^m(a_n)-a_n|_2>0$, for all $m\neq 0$. Notice that for $g\neq h$ in a i.c.c 
group, we have $|g-h|_2=2^{1/2}$, hence it suffices to prove that  $\alpha^m(a_n)$ is not equal to 
$a_n$, for sufficiently large $n$. Using the notation in the proof of Theorem \ref{t1} (see relation 
(\ref{e1})) we get $\epsilon_{a_n}=d_n$. If  $\alpha^m(a_n)=a_n$ then we obtain 
$$\epsilon_{x_0^ma_nx_0^{-m}}=\epsilon_{a_n}$$  Applying (\ref{e1}) $$x_0^m(d_n)=d_n$$ The last equality 
cannot happen though, as $x_0^m$ has slope equal to $N^m$ near $x=1$ and $d_n\rightarrow 1$. 
\end{proof}

\begin{corollary}For any integer $N\geq 2$, the generalized Thompson group $F(N)$ is inner amenable.
\end{corollary}
\begin{proof}From the theorems above $\mathcal{L}(F(N))$ has property $\Gamma$. 
\end{proof}
The next proposition establishes an exact sequence that allows us once again to conclude $F(N)$ is 
inner amenable. For $N=2$ this proposition specializes in Theorem 4.1. of \cite{Can} (except the part 
about the commutator subgroup, see also \cite{Che}) and also appears in a more general form in \cite{St}.     
\begin{proposition}\label{pr2}One has the short exact sequence
$$1\rightarrow F^{'}\rightarrow F(N)\rightarrow\mathbb{Z}^2\rightarrow 1$$
\end{proposition} 
\begin{proof}We prove that the following group morphism is onto: $\phi:F(N)\rightarrow \mathbb{Z}^2$, $\phi(f)=(a,b)$, 
where $f$ has slope $N^a$ near $x=0$ and slope $N^b$ near $x=1$. Suffices to show that there exist $f_1$ and $f_2$ in $F(N)$ 
such that $\phi(f_1)=(1,0)$ and $\phi(f_2)=(0,1)$. Let $p>0$ and $d:=1/N^p$ such that $d(N+1)<1$. Define now 

$$f_1(x)=\left\{
\begin{array}{lr}
Nx,&\mbox{ } 0\leq x\leq d\\
\frac{x}{N}+Nd-\frac{d}{N} ,&\mbox{ } d\leq x\leq d(N+1)\\
x,&\mbox{ } d(N+1)\leq x\leq 1\\
\end{array}\right.$$
 
Clearly $f_1\in F(N)$ and $\phi(f_1)=(0,1)$; $f_2$ can be obtained by applying a symmetry to $f_1$. 
In conclusion $\phi$ is onto. Notice that its kernel is exactly the normal subgroup $F^{'}$. 
\end{proof}
\par
Question: Following \cite{Br} in a particular case we obtain that $F(2)=F$ is not isomorphic to $F(3)$. 
Is it true that $\mathcal{L}(F(2))\cong\mathcal{L}(F(3))$? Notice that if these factors are 
not isomorphic then by the uniqueness of the hyperfinite $II_{1}$ factor at least one of the groups $F$ or $F(3)$ 
follows non amenable. 
\begin{acknowledgements}We would like to thank Dorin Dutkay for many useful suggestions. 
We also thank prof. Sean Cleary for indicating the right references on the generalized Thompson group. 
We express our gratitude to prof. Florin Radulescu for his constant support.
\end{acknowledgements}

\end{document}